\documentclass[12pt]{article}

\usepackage[dvips]{epsfig}
\usepackage{amssymb}
\usepackage{graphicx}
\usepackage[T1]{fontenc}
\usepackage{amsmath}
\usepackage{amsthm}
\usepackage{float}

\usepackage{graphicx,amssymb,amsfonts,latexsym,amsmath,amsthm,times}
\usepackage{epsfig}
\usepackage{color}
\usepackage[english]{babel}

\begin{document}
\newtheorem{theo}{Theorem}[section]
\newtheorem{prop}[theo]{Proposition}
\newtheorem{lem}[theo]{Lemma}
\newtheorem{coro}[theo]{Corollary}
\theoremstyle{definition}
\newtheorem{definition}[theo]{Definition}
\newtheorem{example}[theo]{Example}
\newtheorem{xca}[theo]{Exercise}

\numberwithin{equation}{section}

\newcommand{\norm}[1]{\lVert#1 \rVert}
\newcommand{\be}{\beta(\epsilon,c,y)}
\newcommand{\ba}{\overline{\beta}(\ep,c)}
\newcommand{\na}{\nabla}
\newcommand{\ep}{\epsilon}
\newcommand{\R}{\mathbb{R}}
\newcommand{\La}{\Lambda}
\newcommand{\la}{\lambda}
\newcommand{\h}{\tilde{H}}
\newcommand{\g}{\tilde{G}}
\renewcommand{\k}{\tilde{K}}
\newcommand{\w}{\tilde{W}}
\renewcommand{\eqref}[1]{(\ref{#1})}\

\begin{center}
{\Large {\bf Numerical Simulations of Heat Explosion With Convection In
Porous Media}}
\end{center}

\vspace*{0.5cm}
\begin{center}
Karam Allali$^a$\footnote{Corresponding author. E-mail:
allali@fstm.ac.ma}, Fouad Bikany$^a$, Ahmed Taik$^a$
 and Vitaly Volpert$^b$
\end{center}

\begin{center}
$^{a}$ Department of Mathematics, FSTM, MAC Laboratory, University Hassan II, Po.Box146, Mohammedia, Morocco\\
$^{b}$ Department of Mathematics, ICJ, University Lyon 1, UMR 5208
CNRS, Bd. 11 Novembre, 69622 Villeurbanne, France
\end{center}

{\bf Abstract:} In this paper we study the interaction between
natural convection and heat explosion in porous media. The model
consists of the heat equation with a nonlinear source term
describing heat production due to an exothermic chemical reaction
coupled with the Darcy law. Stationary and oscillating convection
regimes and oscillating heat explosion are observed. The models
with quasi-stationary and unstationary Darcy equation are
compared.

\vspace{.2cm}


\vspace{.2cm}





\section{Introduction}

%

 The theory of heat explosion began from the classical works by
 Semenov \cite{s}  and Frank-Kame\-netskii \cite{fk}.
 In Semenov's theory, the temperature distribution in the vessel is supposed to be
 uniform. An average temperature in the vessel is described by the ordinary differential equation

 \begin{equation}
 \label{int1}
 \frac{d \theta}{dt} =  e^\theta - \lambda \theta .
 \end{equation}
 The first term in the right-hand
 side corresponds to heat production due to an exothermic chemical
 reaction, the second term to heat loss through the
 boundary of the vessel.
 In Frank-Kamenetskii's theory, spatial temperature distribution is
 taken into account. The model consists of the reaction-diffusion
 equation,

 \begin{equation}
 \label{fk}
 \frac{d \theta}{dt} =  \Delta \theta + F_K e^\theta ,
 \end{equation}
 where the first term in the right-hand side describes heat
 diffusion, $F_K$ is called the Frank-Kame\-netskii
  parameter. This equation is considered in a bounded domain with
  the zero boundary condition for the dimensionless temperature.

  In both models, heat explosion was treated as an unbounded
 growth of temperature (blow-up solution). Thus, the problem of
 heat explosion was reduced to investigation of existence,
 stability and bifurcations of stationary solutions of
 differential equations. These questions initiated a big body of
 physical and mathematical literature (see \cite{zblm} and the references
 therein).

 The effect of natural convection on heat explosion was first studied in \cite{ksp,ms}.
 It was shown that the critical value of the Frank-Kamenetskii
  parameter increases with the Rayleigh number and explosion can be prevented by vigorous
  convection.
 These works were continued by \cite{bv,dv,dgmv,lvm}
 where new stationary and oscillating regimes were found. The authors showed how complex regimes appeared through successive bifurcations
   leading from a stable stationary temperature distribution without convection to a stationary symmetric convective solution,
    stationary asymmetric convection, periodic in time oscillations, and finally aperiodic oscillations.
    Oscillating heat explosion, where the temperature grows and
    oscillates, was discovered.
  The effects of natural convection and
  consumption of reactants on heat explosion in a closed spherical vessel were studied in \cite{liu}.
  The influence of stirring on the limit of thermal explosion was investigated in \cite{kbjs}.
 Heat explosion with convection in a horizontal cylinder was
 considered in \cite{our}.

 All these works study heat explosion in a gaseous or liquid
 medium with its motion described by the Navier-Stokes equations
 under the Boussinesq approximation. Thermal ignition in a porous
 medium is investigated in \cite{Kord}. The Darcy law in
 a quasi-stationary form under the Boussinesq approximation is
 used to describe fluid motion. It is shown that
 convection decreases the maximal temperature and increases the
 critical value of the Frank-Kamenetskii parameter.
 The interaction of free convection and exothermic chemical is
 studied in \cite{vh1}. The authors consider zero-order exothermic
 reaction in a rectangular domain and find the onset of convection
 by an approximate analytical method. Similar problem with depletion
 of reactants is investigated in \cite{vh3}. Ignition time of heat
 explosion in a porous medium with convection is found in
 \cite{vh2}. Heat explosion in one-dimensional flow in a porous
 medium is studied in \cite{balak}.

 In this work we continue to study interaction of natural convection with thermal explosion in porous
 media. The reaction-diffusion equation for the temperature distribution will
 be coupled with Darcy's law describing fluid motion in porous
 media. Along with stationary regimes with and without convection,
 we will show the existence of oscillating convective regimes and
 oscillating heat explosion, which were not yet observed for this
 problem. We will also compare two models, with quasi-stationary
 approximation and complete Darcy equation.

The paper is organized as follows. The first problem with a
quasi-stationary Darcy equation is formulated in Section 2,
followed in Section 3 by numerical simulations. Section 4 is
devoted to the second model with the complete Darcy equation.
Short conclusions are given in the last section.


\section{Governing equations}


We consider the first order reaction,

\begin{equation}
A \stackrel{K(T)}{\longrightarrow} B,
\end{equation}
 and the temperature dependence of the reaction rate $K(T)$ given by the
 Arrhenius law:

\begin{equation}
K(T) = k_0 exp\left(-\frac{E}{R T}\right),
\end{equation}
where $E$ is the activation energy, $T$ the temperature, $R$ the universal gas constant and $k_0$ the pre-exponential factor.
 The model consists of the reaction-diffusion equation with
 convective terms and of Darcy's law in the quasi-stationary
 approximation for an incompressible fluid:

\begin{equation}
\frac{\partial T}{\partial t}+\mathbf{v}.\nabla T
 = \kappa \Delta T+ q K(T)
\end{equation}

\begin{equation}
 \label{por1}
\mathbf{v}+\frac{K}{\mu}\nabla p = \frac{g\beta K}{\mu}\rho(T-T_0) \mathbf{\gamma }
\end{equation}

\begin{equation}
 \label{por2}
\nabla .\mathbf{v} = 0.
\end{equation}
Here $\mathbf{v}$ denotes the fluid velocity field, $p$ is the pressure, $\kappa$ the coefficient of thermal diffusivity,
 $\mu$ the kinematic viscosity, $\rho$ the density, $q$ the heat release, $g$ is the acceleration due to gravity,
  $\gamma$ a unit vector in the vertical direction, $T_0$ the characteristic value of the temperature, $K$ the
  permeability. Depletion of reactants in the heat balance equation is
  neglected. It is a conventional assumption in the theory of heat
  explosion.
%
%
%
 This system is considered in a $2D$ square domain, $0 \le x \le 2\ell$, $0 \le y  \le 2\ell$.
 The boundary conditions will be specified below.

 Equations of motion (\ref{por1}), (\ref{por2}) are written under
 the Boussinesq approximation and quasi-stationary approximation.
 The former signifies that the density of the fluid is everywhere
 constant except for the buoyancy term (right-hand side in
 (\ref{por1})). Quasi-stationary approximation in the Darcy law,
 though often used for fluids in porous medium, should be verified
 for the problem of heat explosion. We return to this question in
 Section 4.


In order to write the dimensionless model, we introduce new
spatial variables $\displaystyle x'=x/\ell$, $\displaystyle
y'=y/\ell$, time
 $\displaystyle t'=\frac{\kappa}{\ell^2}t$, velocity $\displaystyle \frac{\ell}{\kappa}\mathbf{v}$ and pressure $\displaystyle \frac{K}{\mu \kappa} p$.
  Denoting $\displaystyle \theta = \frac{E(T-T_0)}{RT_0^2}$ and keeping for convenience the same notation for the other variables, we obtain:

\begin{equation}\label{sd1}
\frac{\partial \theta}{\partial t}+u\frac{\partial \theta}{\partial x}+v\frac{\partial \theta}{\partial y}
 = \frac{\partial^2 \theta}{\partial x^2}+\frac{\partial^2 \theta}{\partial y^2}+F_K e^{\theta}
\end{equation}

\begin{equation}\label{sd11}
u+\frac{\partial p}{\partial x}=0
\end{equation}

\begin{equation}\label{sd22}
v+\frac{\partial p}{\partial y}
= R_p \theta
\end{equation}

\begin{equation}\label{sd2}
\frac{\partial u}{\partial x}+\frac{\partial v}{\partial y}
 = 0.
\end{equation}
Here $\displaystyle F_K = \frac{Eqk_0\alpha e^{-\frac{E}{RT_0}}
\ell^2}{RT_0^2\kappa}$ is the Frank-Kemenetskii parameter,
$\displaystyle R_p = \frac{K\rho R_a}{\ell^2}$, $\displaystyle
R_a=\frac{g\beta R T_0^2 \ell^3}{E\kappa \mu}$ is the Rayleigh
number, $(u,v)$ is the velocity vector. Under the assumptions of
large activation energy, $RT_0/E\ <<\ 1$, we can perform the
Frank-Kamenetskii transform, so that
 the nonlinear reaction rate in the equation \eqref{sd1} is taken to be $F_K exp(\theta)$ \cite{fk,zblm}.

 Let us note that the characteristic thermal diffusion time scale, $\kappa/l^2$,
 which enters the dimensionless variables, is related to the rate
 of heat loss through the boundary. Competition of heat loss
 with heat production due to chemical reaction determines
 conditions of heat explosion. The presence of convection,
 which intensifies heat loss, provides an additional factor that
 can influence heat explosion.

The system of equations \eqref{sd1}-\eqref{sd2} is supplemented by
the boundary conditions:

\begin{equation}
x=0, 2  : \ \ \ \   \frac{\partial \theta}{\partial x} =u=0 ,
\end{equation}

\begin{equation}
y=0, 2  : \ \ \ \ \theta=v= 0 .
\end{equation}


\section{Numerical simulations}

\subsection{Numerical method}

To describe the numerical method, we first introduce the stream
function $\psi$ using incompressibility of the fluid:

\begin{equation}
\displaystyle \left(\begin{array}{ll} u \\ v \end{array}\right) =
\left(\begin{array}{ll}\displaystyle\  \frac{\partial
\psi}{\partial y} \\ \displaystyle  -\frac{\partial \psi}{\partial
x} \end{array}\right) .
\end{equation}
We apply the rotational operator to equations
\eqref{sd11}-\eqref{sd22} in  order to eliminate the pressure. The
equation for the stream function writes

\begin{equation}\label{3}
-\Delta \Psi = R_p \partial _{x}\theta .
\end{equation}
 From the boundary conditions for the velocity we obtain the
 boundary conditions for the stream function:

 $$ x,y=0,2 : \Psi = 0 . $$
 This problem is solved by the fast Fourier transform taking into account the
Dirichlet boundary conditions. Equation \eqref{3} is coupled to
equation \eqref{sd1}. The latter is solved using an implicit
finite difference scheme and alternative direction method. It is a
simple and robust method often used for reaction-diffusion
problems.
%
%
%
%


\subsection{Results}


 If the fluid velocity in the porous medium is zero, $u=v=0$, then
 system of equations (\ref{sd1})-(\ref{sd2}) is reduced to the
 single reaction-diffusion equation (\ref{fk}). If the Frank-Kamenetskii parameter
 is less than the critical value $F_K=1$, then there are two
 stationary solutions which depends only on the
 vertical variable. If $F_K>1$, then stationary solutions do not
 exist, and the solution of the evolution problem grows to
 infinity. This case corresponds to heat explosion.

\begin{figure}[!t]
\centering
 \includegraphics[scale=0.5]{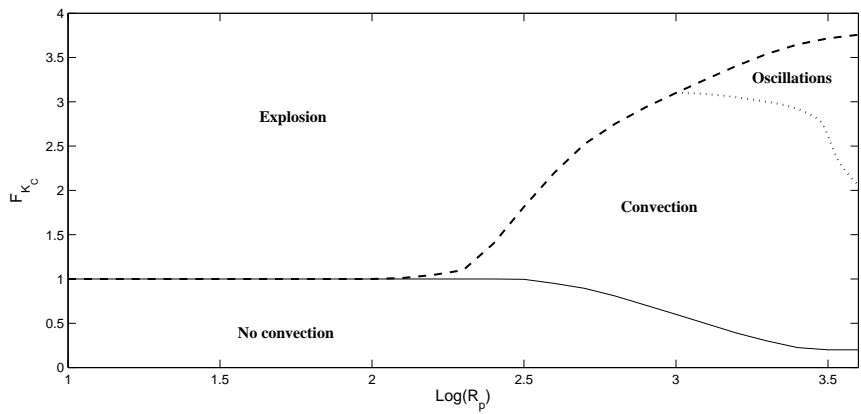}
\caption{Numerical simulations reveal four different regions in
the parameter plane $(R_p,F_K)$. In the lower region, there are
stationary regimes without convection. In the upper region, they
do not exist and solution blows up. In the intermediate region
from the right, the solution remains bounded (no explosion). It
can be stationary or oscillating with convection. }
\end{figure}

 Convection can change conditions of heat explosion. Figure 1
 shows four domains in the plane of parameters $(F_K,R_p)$:
 stationary regimes without convection, bounded (stationary or
 oscillating) solutions with convection, blow-up solutions
 (explosion).

\begin{figure}[!t]
\centering
 \includegraphics[scale=.3]{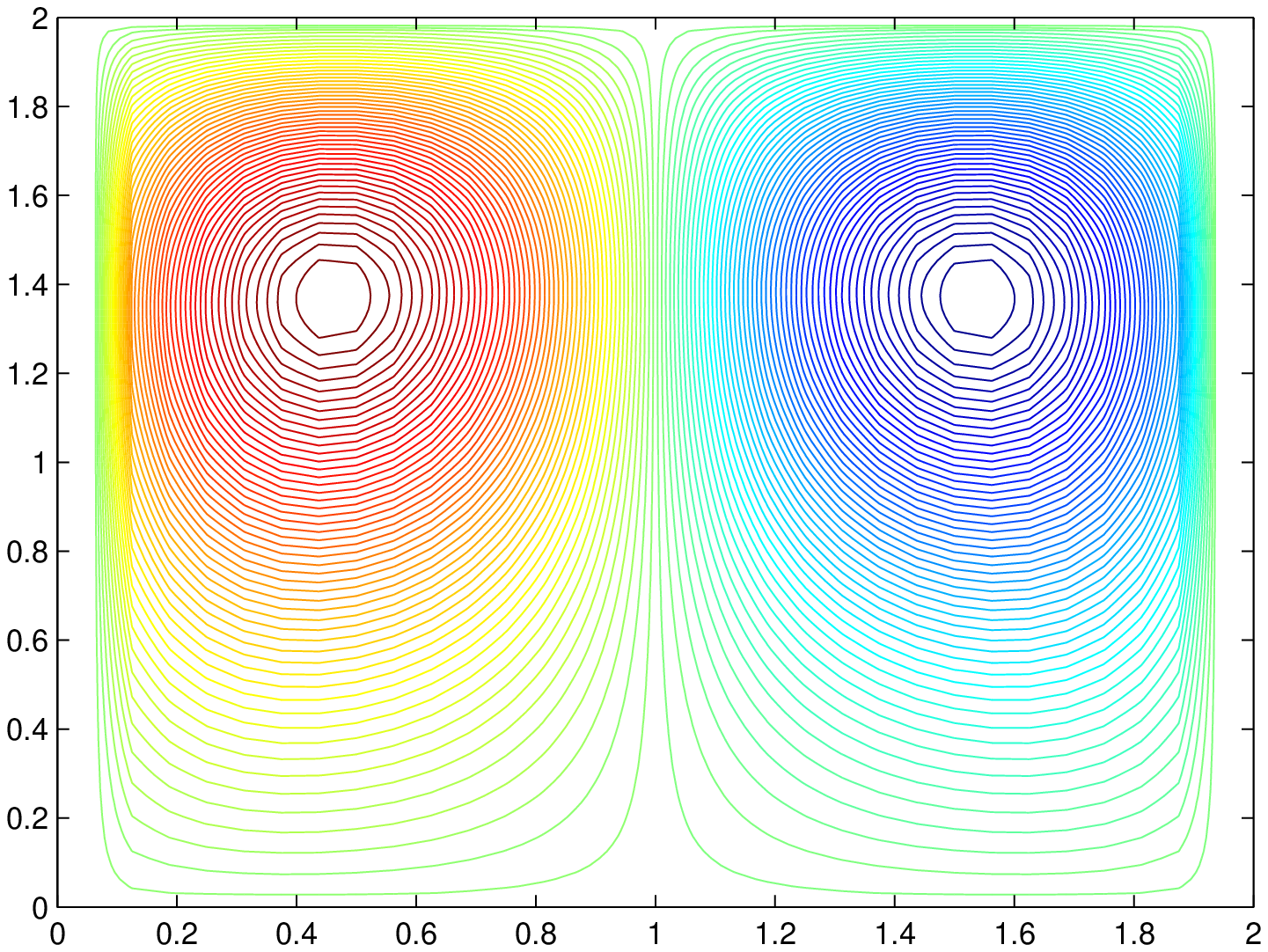} \includegraphics[scale=.3]{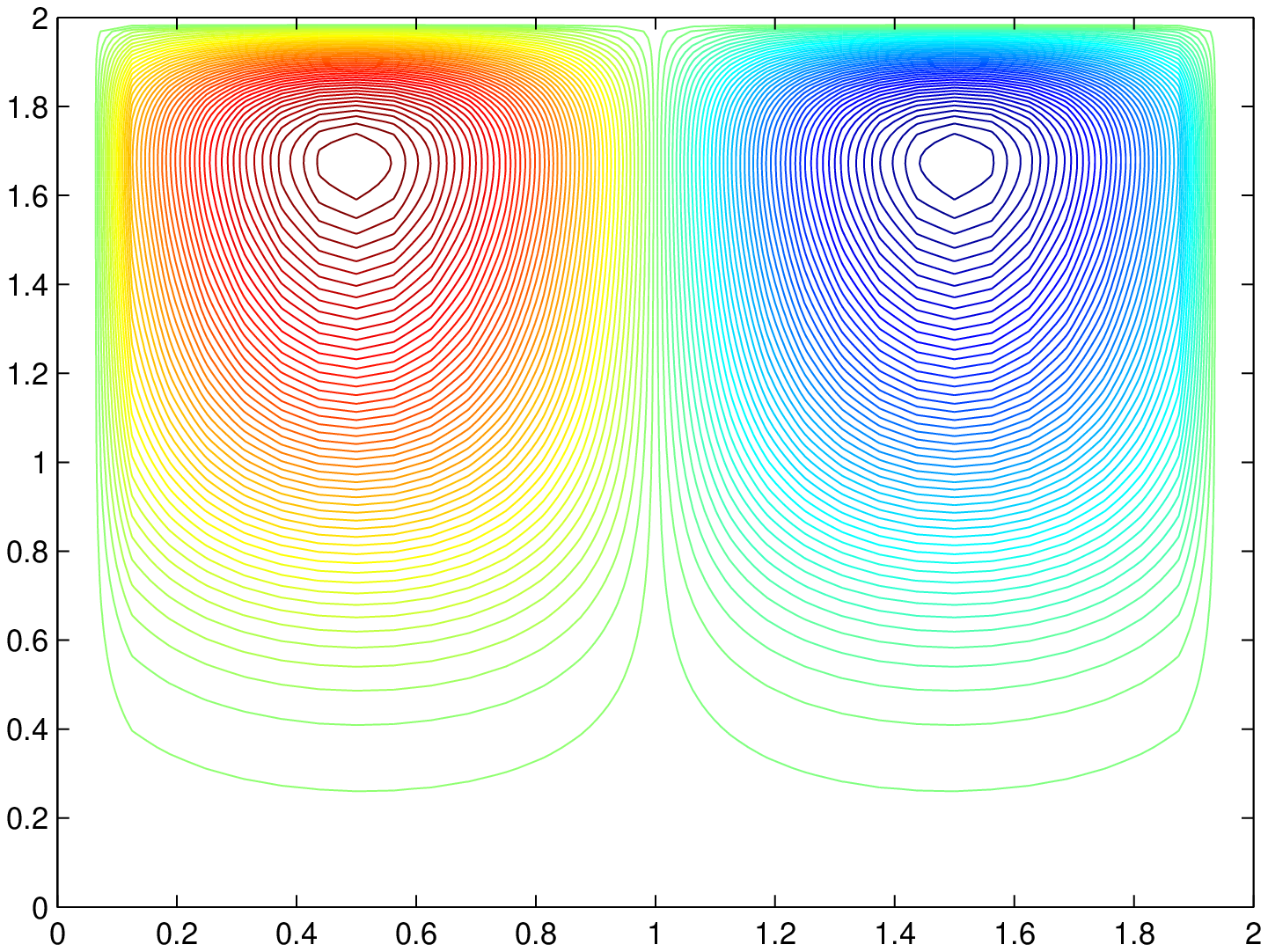} \includegraphics[scale=.3]{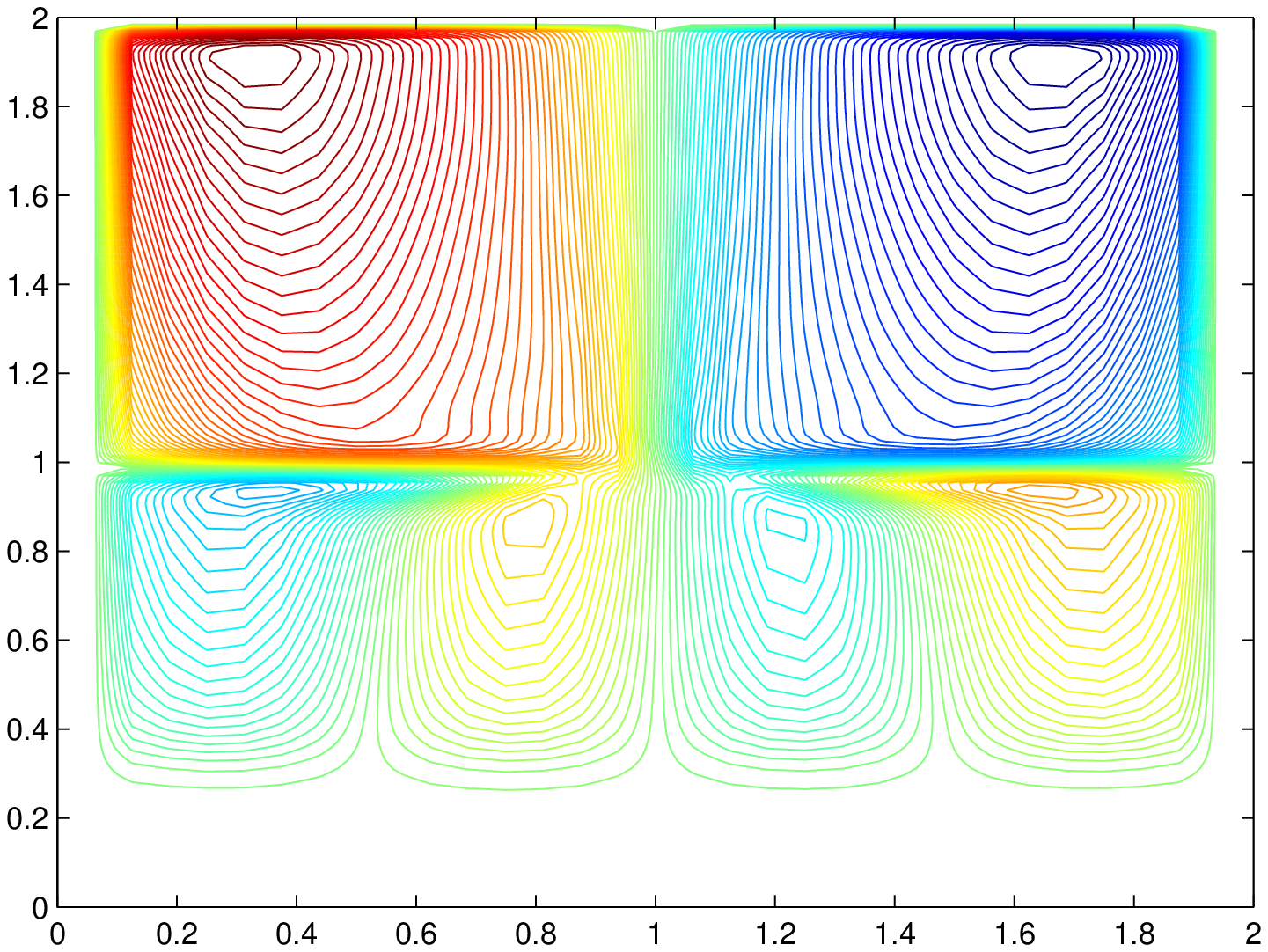}
\caption{Stationary convective regimes. Level lines of the stream
function for $F_K = 1$, $R_p=100$ (left), for $F_K = 1$,
$R_p=1000$ (middle) and for $F_K = 3.2$, $R_p=1000$ (right).}
\end{figure}


\begin{figure}[!t]
\centering
 \includegraphics[scale=.3]{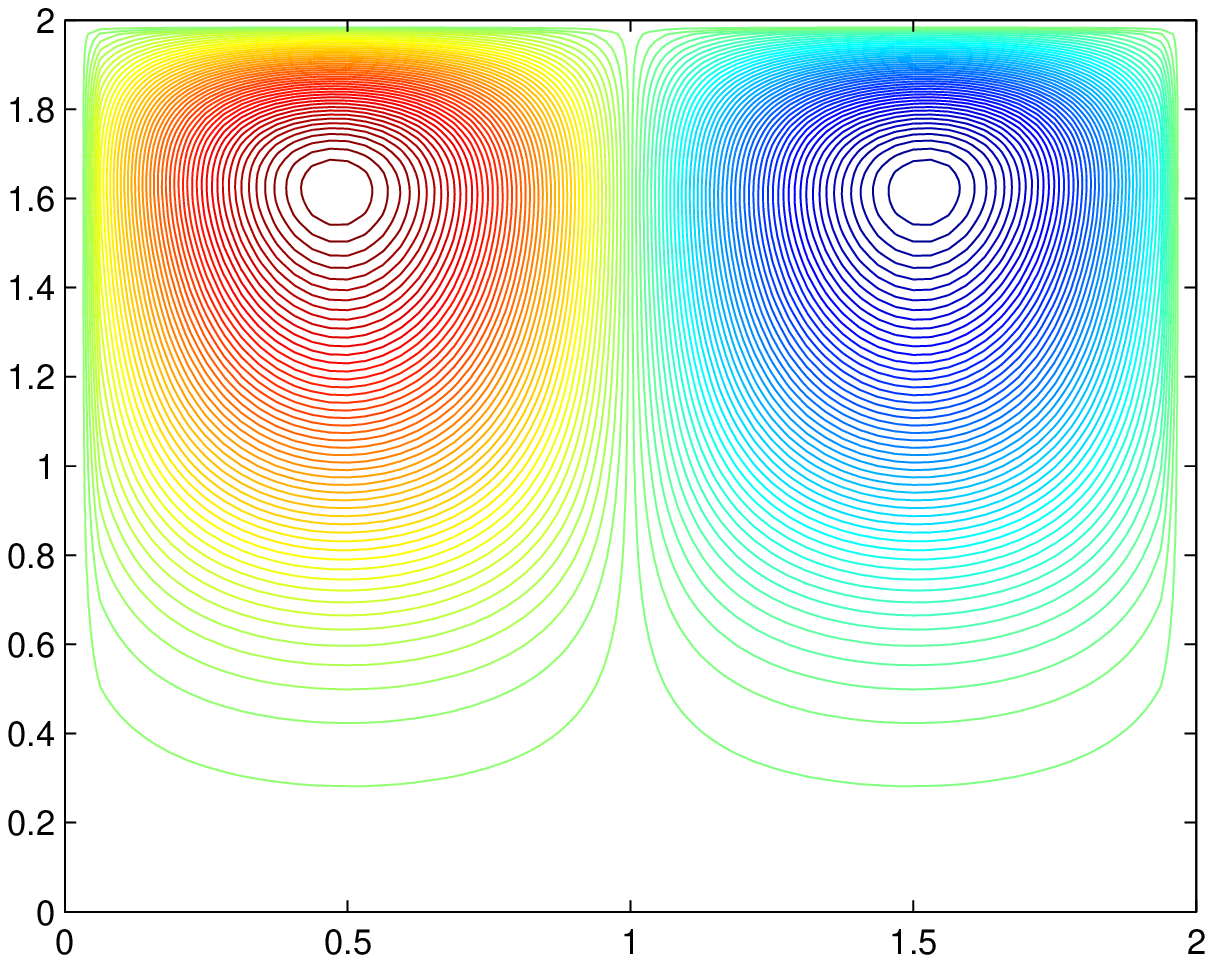} \includegraphics[scale=.3]{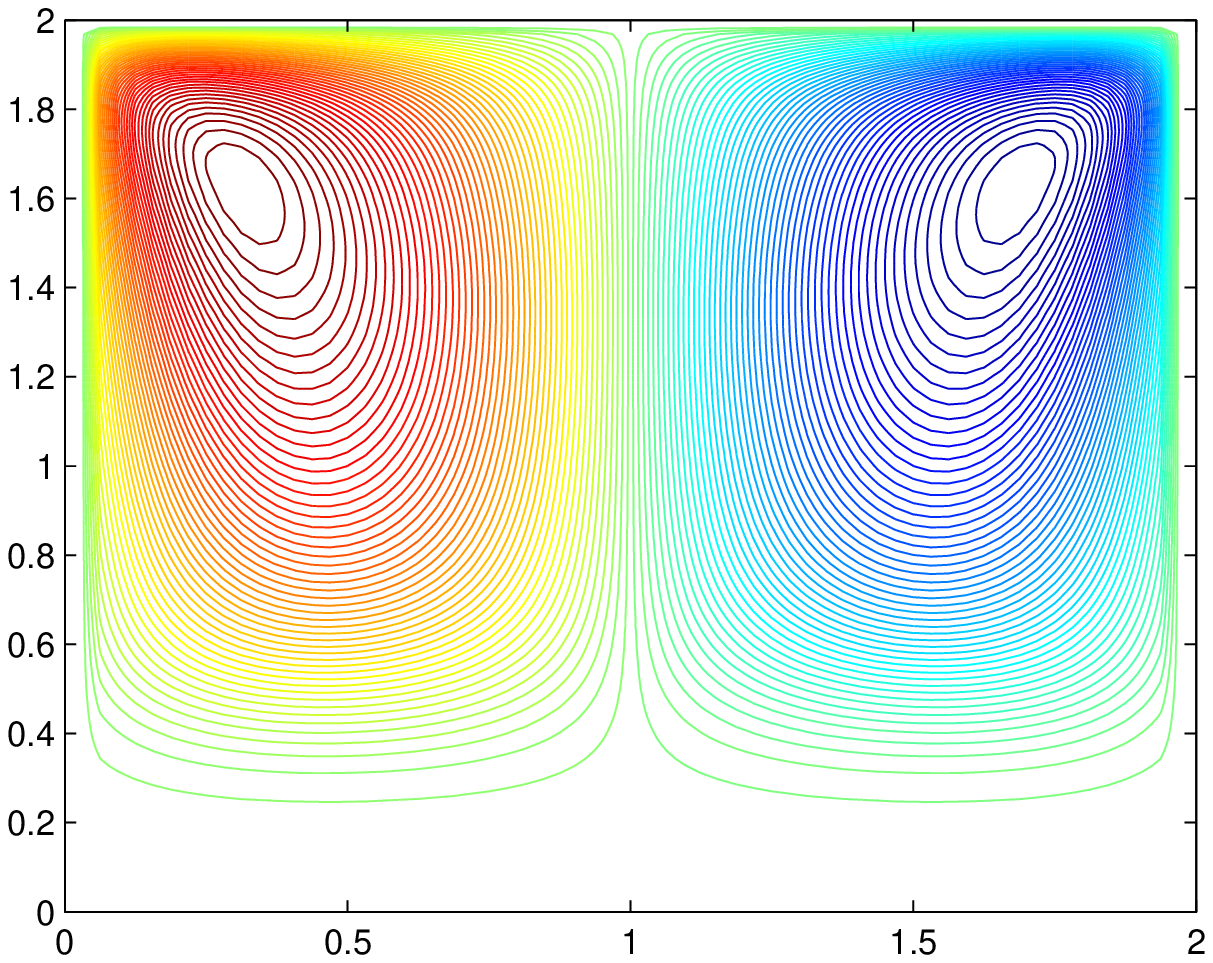}
 \includegraphics[scale=.3]{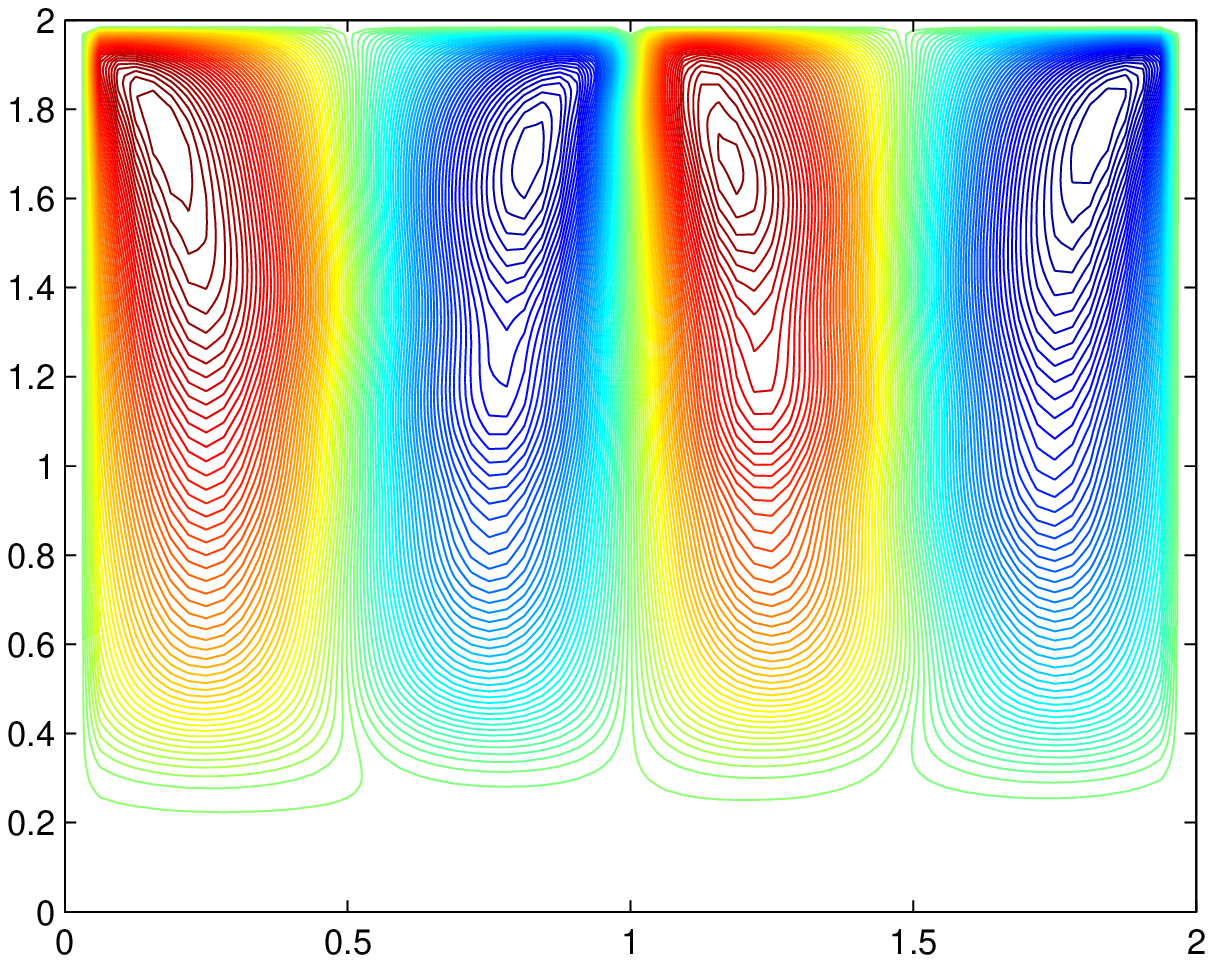}
\caption{Stationary and Oscillating convective regimes. Level
lines of the stream function for $R_p=3500$ and for $F_K = 0.2$
(left), $F_K = 1$ (middle) and for $F_K = 3.75$ (right).}
\end{figure}

 If the Rayleigh number is sufficiently small, then there is no
 convection and the critical value of the Frank-Kamenetskii
 parameter is independent of $R_p$.
 If we fix $F_K$ and increase $R_p$, then the stationary solution
 without convection loses its stability resulting in appearance of
 stationary convective regimes.

  We note that convection increases heat exchange through the boundary.
 Therefore, the critical value of $F_K$, for which explosion occurs, grows when $R_p$
 increases.
 If $F_K$ is sufficiently large, the temperature becomes
 unbounded, which corresponds to heat explosion.
 Qualitatively this diagram is similar to the case of heat
 explosion with convection in fluids \cite{ms,dgmv}.

\begin{figure}[!t]
\centering
$\begin{array}{c}
\includegraphics[scale=.43]{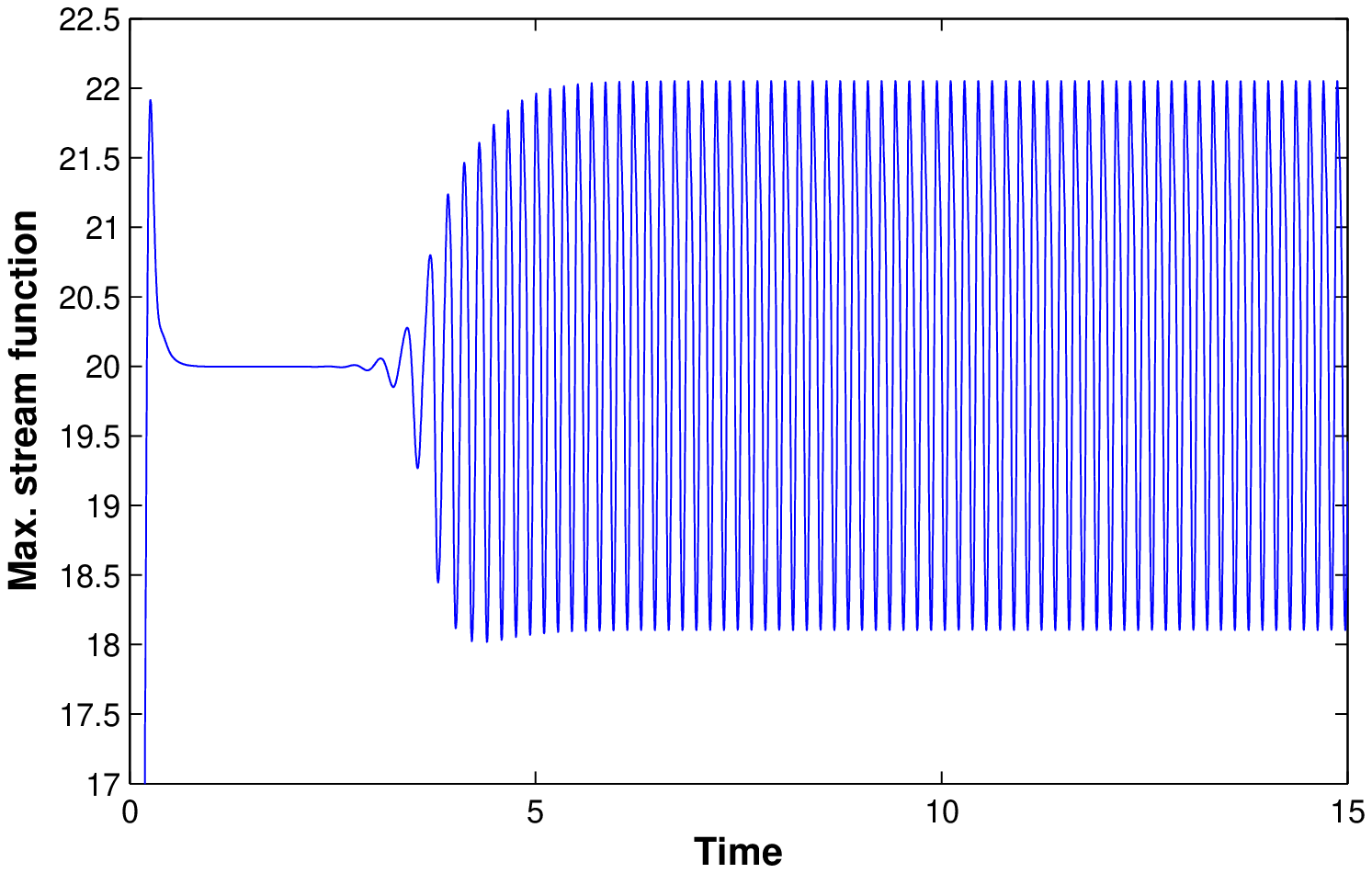}\vspace{-.3cm}\\
(a)\vspace{-.17cm}\end{array}$
 $\begin{array}{c}
 \includegraphics[scale=.35]{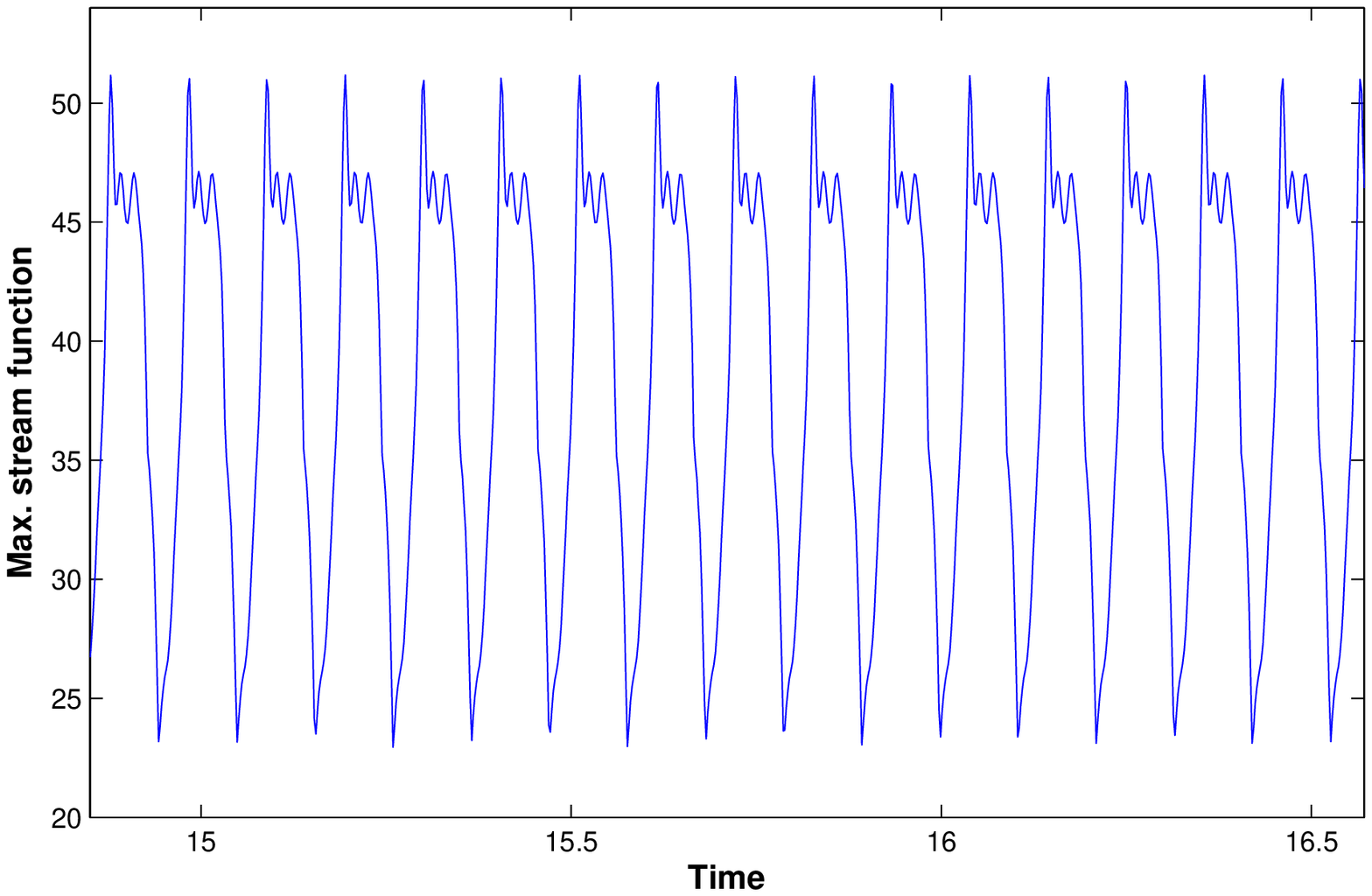}\vspace{-.3cm}\\
(b)\vspace{-.17cm}\end{array}$
 \caption{Oscillating convective regimes. Maximum of the stream
function as a function of time for $R_p=4\times10^3$ and for
different values of $F_K$, $(a) F_K=2.2$, $(b) F_K=3.7$.}
\end{figure}

\begin{figure}[!t]
\centering
$\begin{array}{c} \includegraphics[scale=.37]{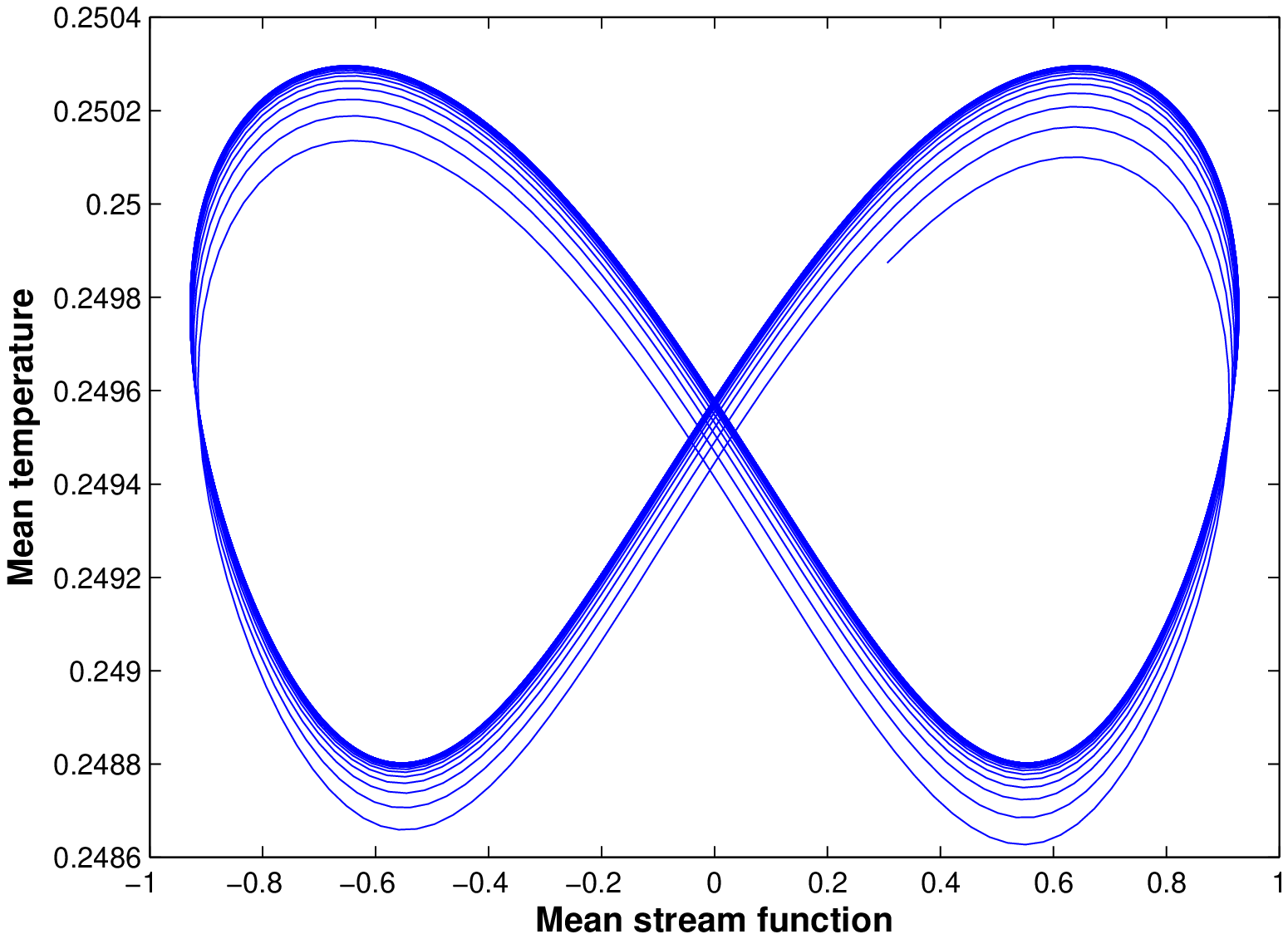}\vspace{-.3cm}\\  (a)\vspace{-.17cm} \end{array}$
$\begin{array}{c}  \includegraphics[scale=.35]{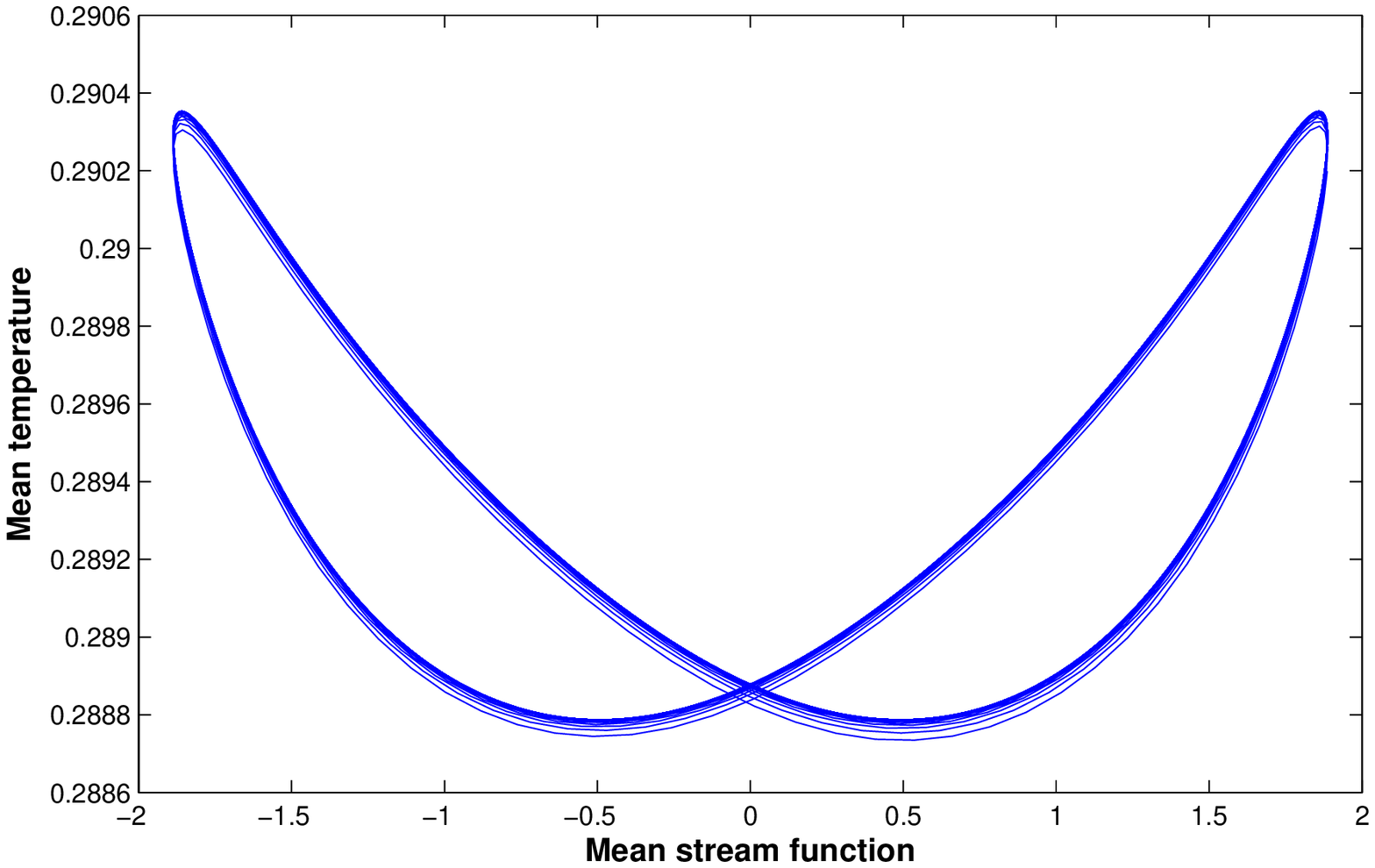}\vspace{-.3cm}\\  (b)\vspace{-.17cm}\end{array}$\\
$\begin{array}{c}
\includegraphics[scale=.35]{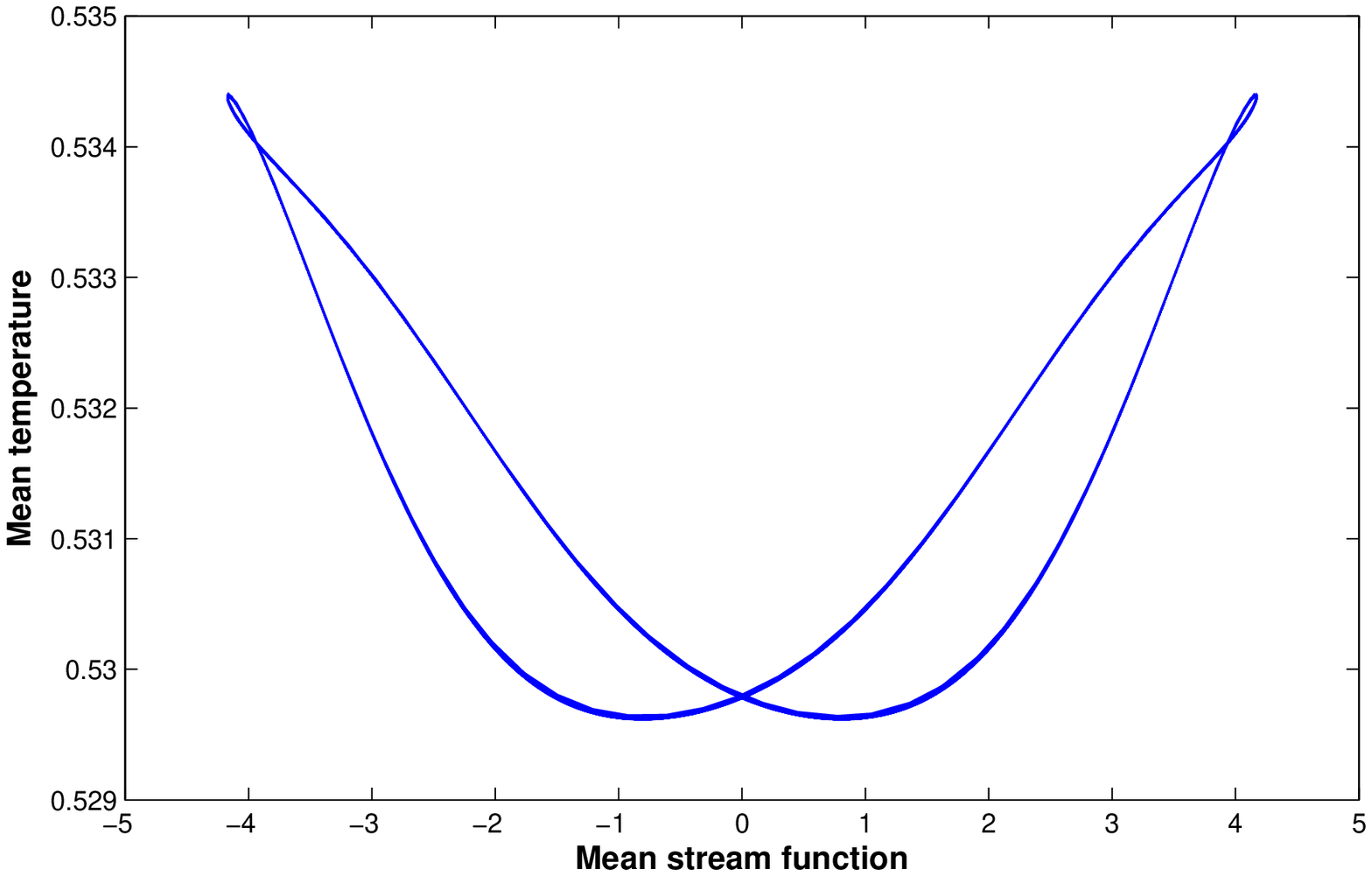}\vspace{-.3cm}\\
(c)\vspace{-.17cm}\end{array}$ $\begin{array}{c}
\includegraphics[scale=.35]{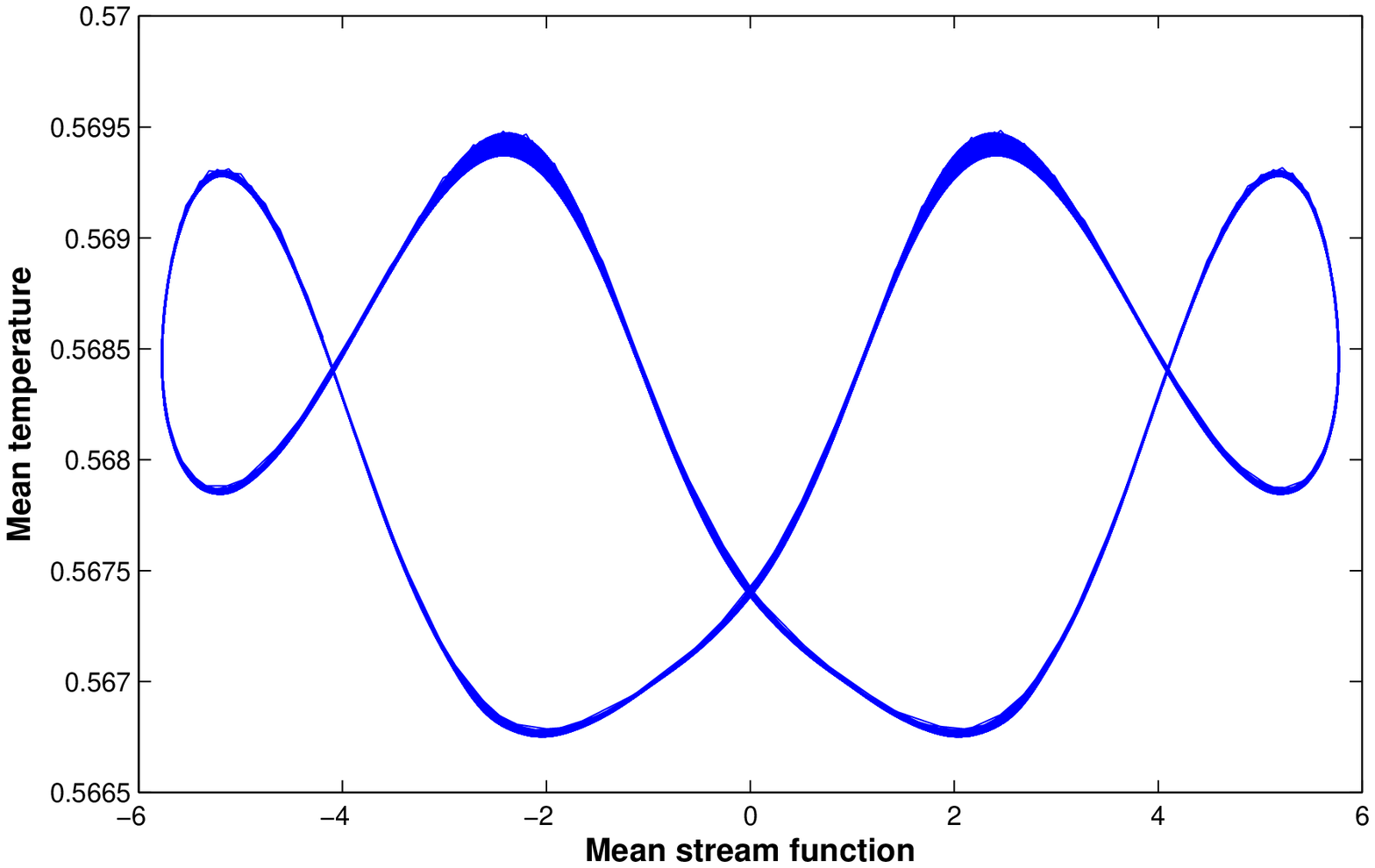}\vspace{-.3cm}\\
(d)\vspace{-.17cm}\end{array}$ \caption{Oscillating convective
regimes.  Mean of temperature as function of the mean value of the
stream function for $R_p=4\times10^3$ and for different values of
$F_K$, $(a) F_K=2$, $(b) F_K=2.2$, $(c) F_K=3.5$, $(d) F_K=3.7$.}
\end{figure}

 Figure $2$ shows level lines of the stream function for
 different values of parameters in the case of stationary
 convection.
 The maximum of the stream function grows significantly with the increase of the
 Rayleigh number. It is about $2\times 10^{-3}$ for $R_p=100$ and about $4.5$ for $R_p=1000$.
 Besides, in the second case the vortices are located closer to the upper boundary
 (Figure 2, left and middle). For higher values of $F_K$ (Figure
 2, right), the maximal temperature increases, convection becomes
 more vigourous and there are additional four vortices located
 below the two main ones.

 Convective patterns depend on the values of parameters.
 Figures $3$ show level lines of the stream function for
 different values of $R_p$ and $F_K$. If we put $R_p=3500$ and
 chose small Frank-Kamenetskii parameter, then we first observe
 two vortex regimes similar to those shown in Figure $2$. The center of the two vortex move to the corners by increasing the
Frank-Kamenetskii parameter (the two first cases correspond to the
stationary convective regime). For larger values of $F_K = 3.75$,
more vortices appear. They compress each other and pack
vertically. This case correspond to the oscillatory convective
regime.



  Figure $4$ shows the maximum of the stream function as a function of time for
  $R_p=4\times10^3$ and for different values of the Frank-Kamenetskii parameter.
    For  $F_k > 2$, instead of a stationary solution, we observe periodic oscillations.
     Figure $5$ shows the mean value of the temperature as a function of the mean
    value of the stream function. The solution forms closed
    curves which correspond to periodic oscillations. The
    structure of these curves changes with the increase of $F_K$.
    For small values, it is a simple $\infty$-shape curve. For
    large values, it becomes double $\infty$-shape. The transition
    from one to the other shown in Figure 5 b), c) where
    additional loops appear at upper corners. More complex
    structure of these curves for large $F_K$ canbe related to
    more complex convective regimes with more vortices.

%

 If we increase the Frank-Kamenetskii parameter even more and cross
 the boundary separating convection and explosion domains in Figure 1,
 the maximum of the temperature and of the stream function oscillate during some time and
 then begin unlimited growth (Figure $6$). It is oscillating heat explosion
 similar to that found before in the case where fluid
 motion was described by the Navier-Stokes equations \cite{dgmv}.
 Thus, along with usual heat explosion where temperature
 monotonically increases, there exists oscillating heat explosion
 where temperature oscillates before explosion. This effect is due
 to the interaction of heat release and natural convection.

\begin{figure}[!t]
 \centerline{\includegraphics[scale=.4]{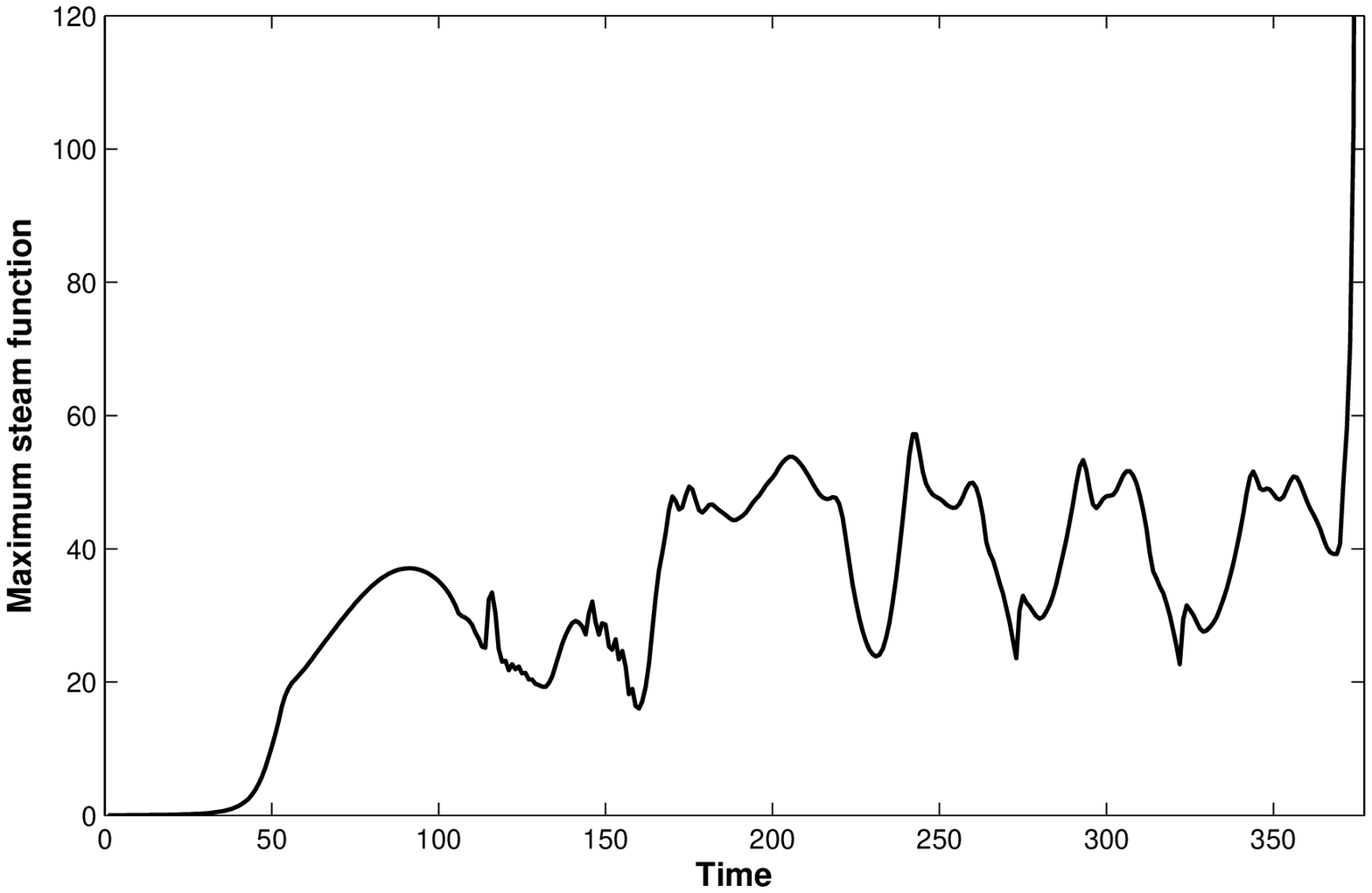}
 \includegraphics[scale=.4]{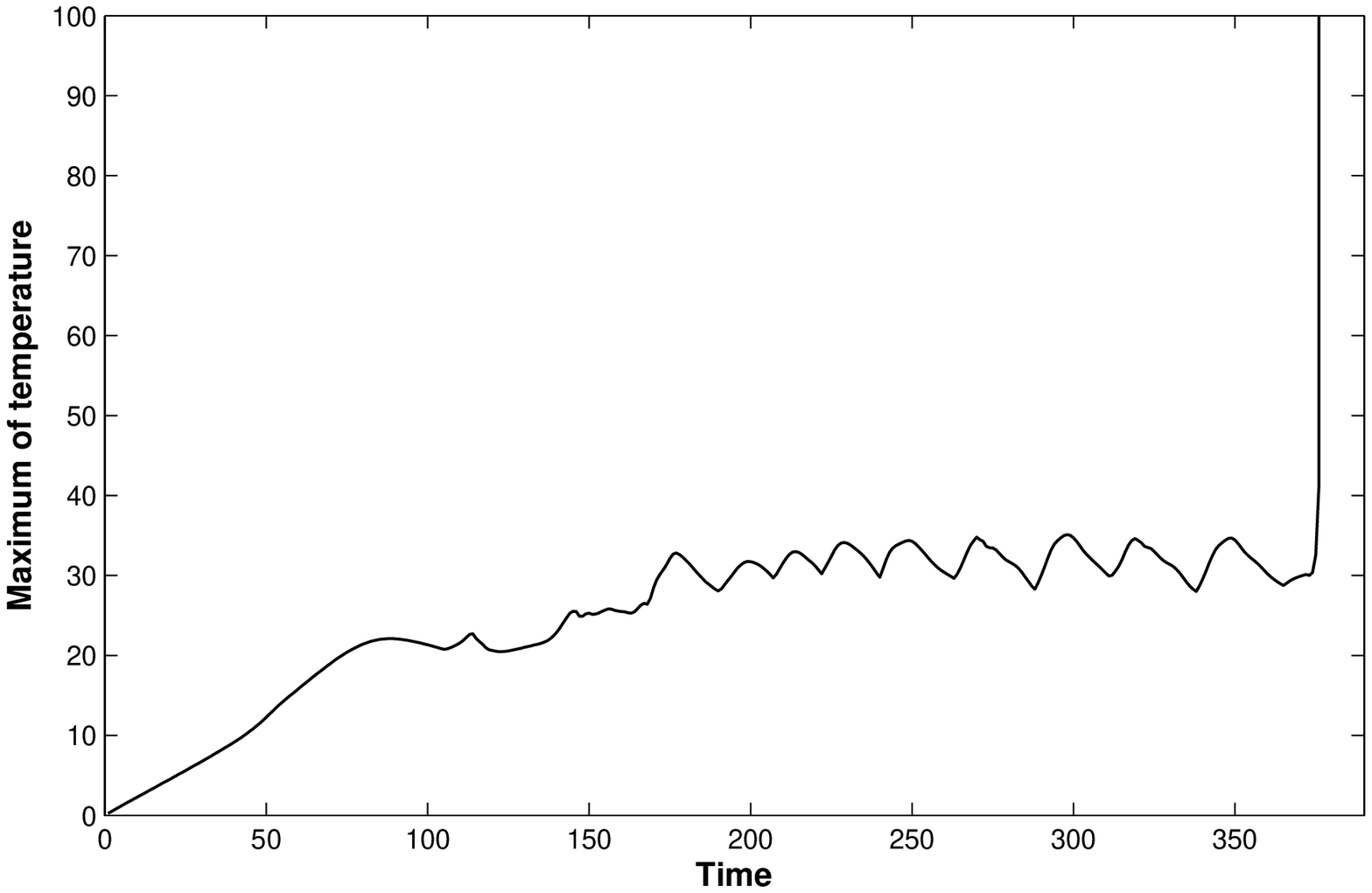}}
\caption{Maximum of the stream function (left) and of the
temperature (right) as functions of time for $R_p=4\times10^3$ and
for $F_K = 3.757$.}
\end{figure}



\section{The model with non-stationary Darcy equation}

\subsection{The new model setting}

In this section, we will consider the following model:

\begin{equation}
\frac{\partial T}{\partial t}+\mathbf{v}.\nabla T
 = \kappa \Delta T+ q K(T)\phi(\alpha),
\end{equation}

\begin{equation}\label{nsd}
 \frac{\partial \mathbf{v}}{\partial t}+\frac{\mu}{K}\mathbf{v}+\nabla p = g\beta\rho(T-T_0) \mathbf{\gamma },
\end{equation}

\begin{equation}
\nabla .\mathbf{v} = 0.
\end{equation}
 The difference in comparison with the model studied above is that
 we do not consider here the quasi-stationary approximation for
 Darcy's law. Introducing vorticity

\begin{equation}
\omega =\it{curl}\ \mathbf{v},
\end{equation}
multiplying the equation (\ref{nsd}) by $\displaystyle
\frac{K}{\mu}$ and following the same steps as in Sections 2 and
3, we obtain the system of equations:

\begin{equation}\label{vsf1}
\frac{\partial \theta}{\partial t}+\frac{\partial \psi}{\partial y}\frac{\partial \theta}{\partial x}-\frac{\partial \psi}{\partial x}\frac{\partial \theta}{\partial y}
 = \frac{\partial^2 \theta}{\partial x^2}+\frac{\partial^2 \theta}{\partial y^2}+F_K e^{\theta},
\end{equation}

\begin{equation}\label{vsf2}
\sigma \frac{\partial \omega}{\partial t}+\omega=R_p\frac{\partial \theta}{\partial x},
\end{equation}

\begin{equation}\label{vsf3}
\omega = -\Delta \psi ,
\end{equation}
 where the parameter $\displaystyle \sigma=\frac{1}{V_a}$ stands for
the inverse of Vadasz number, $V_a = P_r/D_a$, $\displaystyle
P_r=\frac{\mu}{\kappa}$ is the Prandtl number and $\displaystyle
D_a = \frac{K}{\ell^2}$ is the Darcy number.


\subsection{Numerical method and results}

 As before, we use the alternative direction finite difference method to solve equation
 (\ref{vsf1}) and the fast Fourier transform to solve equation (\ref{vsf3}).
 Equation (\ref{vsf2}) is solved by an explicit Euler method.

Figure 7 shows the mean value of the temperature as a function of
the mean value of the stream function for $R_p=4\times10^3$ and
for different values of $\sigma$. When using the first model
(Section 4), we obtain a small "butterfly" in the bottom of the
box (cf. Figure $5$ (a)). This model is a particular case of the
model introduced in Section 4.1 if we formally put $\sigma=0$.
Numerical results in Figure 7 show the convergence of solutions of
system (\ref{vsf1})-(\ref{vsf3}) to the solution of system
(\ref{sd1})-(\ref{sd2}) as $\sigma \to 0$.




\begin{figure}[!t]
\centering
 \includegraphics[scale=0.6]{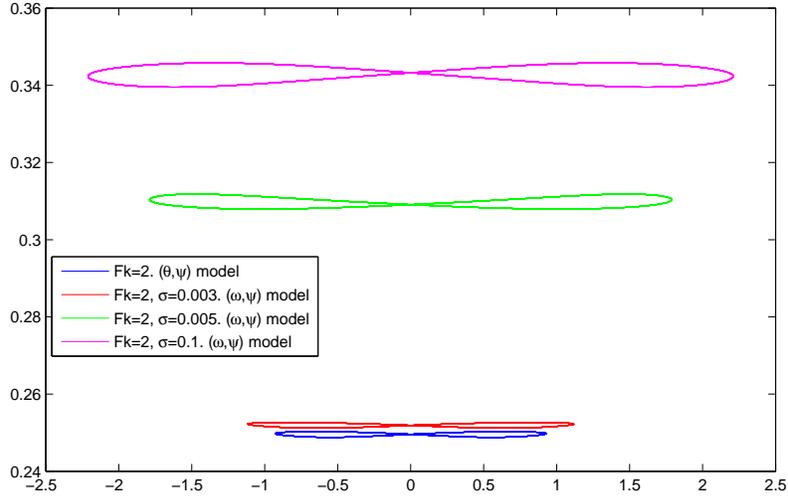}
\caption{Mean of temperature as function of mean of stream
function for $R_p=4\times10^3$ and for different values of
$\sigma$. Solution of the complete Darcy's law converges to the
solution under quasi-stationary approximation as $\sigma$
decreases.}
\end{figure}

\section{Discussion}

In this work, we study the influence of convection on thermal
explosion in a porous medium. We begin with the model where a
nonlinear reaction-diffusion equation is coupled to Darcy's law
written under the quasi-stationary approximation.  Stationary and
oscillating convective regimes are observed. Conditions of
explosion are determined and oscillating heat explosion is found.
We next consider the complete equations of motion without
quasi-stationary approximation. Numerical simulations show the
convergence of the solution to the solution under the
quasi-stationary approximation as the Vadasz number increases.

 Let us note that stationary convective regimes in the problem of
 heat explosion in a porous medium were already observed in
 literature \cite{Kord}, \cite{vh1}. We show in this work that the
 structure of convective solutions depend on the intensity of
 heat release. When the Frank-Kamenetskii parameter increases,
 a second range of vortices appears in the lower part of the
 domain.

 It is interesting that interaction of heat release with
 convection can result in oscillating convective regimes. Some
 indication to this was done in \cite{vh1} where linear stability
 analysis showed that Hopf bifurcation could occur. The
 authors interpreted it as a possible transition to heat explosion
 and not to oscillating convective regimes. On the other hand, we
 have also found oscillating heat explosion where the temperature
 oscillates before it starts unlimited growth. Such regimes can be
 observed if we increase $F_K$ starting from periodically oscillating
 convective solutions.

 Finally, we analyzed applicability of quasi-stationary Darcy law
 used in all previous works devoted to heat explosion in a porous
 medium. It appears that it is a good approximation for large
 values of the Vadasz number. However, it may not be valid for
 certain parameter ranges. In particular, we can expect that
 in some cases the quasi-stationary
 approximation modifies conditions of heat explosion and ignition
 time.

%
%



%
\end{document}